\documentstyle[12pt]{article}
\newcommand{\too}{\longrightarrow}
\newcommand{\om}{\omega}
\newcommand{\Om}{\Omega}
\newcommand{\na}{\nabla}

\newcommand{\al}{\alpha}
\newcommand{\be}{\beta}
\newcommand{\ga}{\gamma}
\newcommand{\Ga}{\Gamma}

\newcommand{\la}{\lambda}

\newcommand{\entraine}{\Longrightarrow}
\newcommand{\inj}{\hookrightarrow}

\newcommand{\ssi}{\Longleftrightarrow}
\def \reel{ {\rm I}\!{\rm R} }

 \def \rat{ {\rm Q}\kern-.65em {}^{{}_/ }}

\newtheorem{Def}{Definition}[section]
\newtheorem{th}{Theorem}[section]
\newtheorem{pr}{Proposition}[section]
\newtheorem{Le}{Lemma}[section]

\title{ Riemann Poisson manifolds and K\"ahler-Riemann foliations }
\author{Mohamed Boucetta} \date{ }\parindent=0cm \begin{document}
\maketitle

{\bf Abstract.} 
\footnote[1]{The authors' research was
supported by The Third World Academy of Sciences RGA No 01-301 RG/Maths/AC}

 Riemann Poisson manifolds were introduced by the author in [1]
and studied in more details in [2].  K\"ahler-Riemann foliations form an
interesting subset of the Riemannian foliations with remarkable properties ( see
[3]).

In this paper we will show that to give a
regular Riemann Poisson structure on a manifold $M$ is equivalent to to give
a K\"ahler-Riemann foliation on $M$  such that the leafwise symplectic form is
invariant with respect to all local foliate perpendicular  vector fields. 
  We show also that the
sum of the vector space of leafwise cohomology and the vector space of the basic
 forms is a subspace of the space of Poisson cohomology.\bigskip

{\it Key words.}   Poisson manifold, Riemannian foliation.

{\it 2000 M. S.  C. Primary 53D17; Secondary 53C12.}

\section{Regular Riemann Poisson manifolds}

 Many fundamental definitions and results about Poisson manifolds can be found in
Vaisman's monograph [7].

Let $P$ be a Poisson manifold with Poisson tensor $\pi$. 
We  have  a bundle map $\pi:T^*P\too TP$ defined by
$$\be(\pi(\al))=\pi(\al,\be),\qquad\al,\be\in T^*P.\eqno(1)$$
On the space of differential 1-forms $\Om^1(P)$, the Poisson tensor induces a Lie
bracket
\begin{eqnarray*}
[\al,\be]_{\pi}&=&L_{\pi(\al)}\be-L_{\pi(\be)}\al-d(\pi(\al,\be))\\
&=&i_{\pi(\al)}d\be-i_{\pi(\be)}d\al+d(\pi(\al,\be)).\qquad\qquad\qquad\qquad\qquad\quad
(2) \end{eqnarray*}

For this Lie bracket and the usual Lie bracket on vector fields, the bundle map
$\pi$ induces a Lie algebra homomorphism $\pi:\Om^1(P)\too{\cal X}(P)$:
$$\pi([\al,\be]_{\pi})=[\pi(\al),\pi(\be)].\eqno(3)$$

Let $<,>$ be a Riemannian metric on the contangent bundle $T^*P$. The
Levi-Civita contravariant connection $D$  associated with the couple $(\pi,<,>)$
is defined by
\begin{eqnarray*} 2<D_\al\be,\ga>&=&\pi(\al).<\be
,\ga>+\pi(\be).<\al ,\ga>-\pi(\ga).< \al,\be>\\
&&+<[\al,\be]_{\pi},\ga >+<[\ga,\al]_{\pi}, \be>+<[\ga,\be]_{\pi}, \al>\qquad(4)
\end{eqnarray*}where $\al,\be,\ga\in\Om^1(P)$.

The notion of contravariant connection was introduced by Vaisman in [7] and
studied in more details by Fernandes in [4].

 The
Levi-Civita contravariant connection $D$  associated with the couple
$(\pi,<,>)$ satisifies:

 $$D_\al\be-D_\be\al=[\al,\be]_{\pi};\eqno(5)$$

 $$\pi(\al).<\be,\ga>=<D_\al\be,\ga>+<\be,D_\al\ga>.\eqno(6)$$

\begin{Def} With the notations above, the triple $(P,\pi,<,>)$ is called 
 a Riemann Poisson manifold if, for any $\al,\be,\ga\in\Om^1(P)$,
$$D\pi(\al,\be,\ga):=
\pi(\al).\pi(\be,\ga)-\pi(D_\al\be,\ga)-\pi(\be,D_\al\ga)=0.\eqno(7)$$

A regular Riemann Poisson manifold is a Riemann Poisson manifold with regular
symplectic foliation. \end{Def}

  Since $D$ has
vanishing torsion and since the contravariant exterior differential $d_\pi$
associated with the bracket $[\;,\;]_\pi$ is given by $d_\pi =-[\pi,.]_S$, we
can deduce obviously that, for any $\al,\be,\ga\in\Om^1(P)$,
 $$-[\pi,\pi]_S(\al,\be,\ga)=D\pi(\al ,\be ,\ga
)+D\pi(\be ,\ga , \al)+ D\pi(\ga , \al, \be).\eqno(8)$$
$[\;,\;]_S$ denote the Shouten bracket.

Let $(P,\pi,<,>)$ be a regular Riemann Poisson manifold.  We denote $TS$ the
involutif distribution tangent to the symplectic foliation.

We have
$$T^*P=Ker\pi\oplus (Ker\pi)^\perp\eqno(9)$$
where $(Ker\pi)^\perp$is the orthoganal of $(Ker\pi)$ with respect to $<,>$.

The Riemannian metric
gives an identification between the cotangent bundle $T^*P$ and the tangent
bundle $TP$ which we denote $\#:T^*P\too TP$.  We put
${\cal H}=\#(Ker\pi).$ We get
$$TP=TS\oplus{\cal H}\quad\mbox{and}\quad \#((Ker\pi)^\perp)=TS.\eqno(10)$$

The bundle
map $\pi$ induces an isomorphism $\pi:(Ker\pi)^\perp\too TS$; we denote 
$\pi^{-1}:TS\too(Ker\pi)^\perp$ its inverse.
The leafwise symplectic form $\om$ belongs
to $\Ga(\wedge^2 T^*S)$ and is given by
$$\om(u,v)=\pi(\pi^{-1}(u),\pi^{-1}(v)),\quad u,v\in TS.\eqno(11)$$

The following proposition gives some elementary properties of the Levi-Civita
contravariant connection.
 \begin{pr} Let $(P,\pi,<,>)$ be a regular Riemann Poisson manifold. Let $D$ be
the  Levi-Civita
contravariant connection associated with $(\pi,<,>)$. Then

1) $\pi(\be)=0$ $\entraine$
$\forall\al\in\Om^1(P),\quad\pi(D_\al\be)=0$;

2) $\pi(\al)=0$ $\entraine$ $D_\al=0;$

3) If $\al,\be\in\Ga((Ker\pi)^\perp)$ then $D_\al\be\in\Ga((Ker\pi)^\perp)$ and
$[\al,\be]_\pi\in\Ga((Ker\pi)^\perp)$.
\end{pr}
{\bf Preuve:} 1) Let $\al,\be,\ga\in\Om^1(P)$ such that $\pi(\be)=0$. We have
$$\ga[\pi(D_\al\be)]=\pi(D_\al\be,\ga)=\pi(\al).\pi(\be,\ga)-D_\al\ga[\pi(\be)]=0.$$

2)  Let $\al,\be\in\Om^1(P)$ such that $\pi(\al)=0$. We have
$D_\al\be=[\al,\be]_\pi+D_\be\al.$
By 1) and (3) it follows that $\pi(D_\al\be)=0.$ Now we show that $D_\al\be\in
\Ga((Ker\pi)^\perp)$. Let $\ga\in\Om^1(P)$ such that $\pi(\ga)=0$.
$\be=\be_1+\be_2$ with $\pi(\be_1)=0$ and $\be_2\in\Ga((Ker\pi)^\perp)$. 

 From (4) and the fact that if $\pi(\al)=0$ and
$\pi(\be)=0$ then $[\al,\be]_\pi=0$, we deduces that $<D_\al\be_1,\ga>=0$.

Finally, $<D_\al\be_2,\ga>=-<\be_2,D_\al\ga>=0$ from (6) and 1).

3) Let $\al,\be\in\Ga((Ker\pi)^\perp)$ and let $\ga\in\Om^1(P)$ 
such that $\pi(\ga)=0$. By (6) and 2) we have
$<D_\al\be,\ga>=-<\be,D_\al\ga>=0$. Hence 
$D_\al\be\in\Ga((Ker\pi)^\perp)$ and since $[\al,\be]_\pi=D_\al\be-D_\be\al$, we
have also $[\al,\be]_\pi\in\Ga((Ker\pi)^\perp)$.$\Box$

\begin{pr} Let $(P,\pi,<,>)$ be a regular Riemann Poisson manifold. We have
$$L_X\pi(\al,\be)=0, \eqno(12)$${for any} $\al,\be\in\Ga((Ker\pi)^\perp)$ { and
for any } $X$ {tangent to} ${\cal H}$.\end{pr}

{\bf Preuve:} For any $\al,\be\in\Ga((Ker\pi)^\perp)$ { and
for any } $X$ {tangent to} ${\cal H}$, we have $[\al,\be]_\pi\in
\Ga((Ker\pi)^\perp)$ and then
$[\al,\be]_\pi(X)=0.$
 A straightforward calculation
 gives the relation 
$$[\al,\be]_\pi(X)=L_X\pi(\al,\be)\eqno(13)$$ and the proposition
follows.$\Box$

We recall the definition of a foliate vector field. A vector field $X\in{\cal
X}(P)$ is said to be foliate if, for all $Y$ tangent to $TS$, $[X,Y]$ is tangent
to $TS$. This is equivalent to the fact that for any $\al\in\Om^1(P)$,
$[X,\pi(\al)]$ is tangent to $TS$.

\section{The correspondence between regular Riemann Poisson manifolds and
K\"ahler-Riemann foliations}
We recall the definition of a K\"ahler-Riemann foliation ( see [3]). A foliation
$\cal F$ on a manifold is called a K\"ahler foliation if it is endowed with a
complex structure $J$ and hermitien metric $h=S-2i\om$ on $T{\cal F}$ such that
$d_{\cal F}\om=0$. A K\"ahler foliation which is also a Riemannian foliation is
called  K\"ahler-Riemann foliation.

 Let $(P,\pi,<,>)$ be a regular Riemann Poisson manifold.

We define a Riemann metric on TP by
\begin{eqnarray*}
g(\#(\al),\#(\be))&=&<\al,\be>,\quad\al,\be\in Ker\pi;\\
g(u,v)&=&<\pi^{-1}(u),\pi^{-1}(v)>,\quad u,v\in TS;\\
g(u,\#(\al))&=&0,\quad \al\in Ker\pi,u\in TS.\end{eqnarray*}

Let $S$ be a symplectic leave. We denote $g_S$ and $\om_S$ the restrictions of
$g$ and $\om$ to $S$. The Levi-Civita connection $\na^S$ of $g_S$ is given by
$$\na^S_{\pi(\al)}\pi(\be)=\pi(D_\al\be),\quad
\al,\be\in\Ga((Ker\pi)^\perp)\eqno(14)$$and we have
$$\na^S\om_S=0.\eqno(15)$$
Furthermore, $S$ is a K\"ahler manifold (see [2]).

We denote $\Om_b^0(P)$ the space of  Casimir functions and
$\Om_b^1(P)$ the space of   basic differential 1-forms. $\al\in\Om_b^1(P)$ if
and only if
$$\pi(\al)=0\quad\mbox{and}\quad i_{\pi(\be)}d\al=0,\quad\forall\be\in\Om^1(P).$$

We have from (2) and a careful verification that
$$\al\in\Om_b^1(P)\quad\ssi\quad\forall\be\in\Om^1(P),\quad
[\al,\be]_\pi=0.\eqno(16)$$

The following proposition gives an interesting characterization of the foliate
vector fields which are tangent to ${\cal H}$.

\begin{pr}  Let $(P,\pi,<,>)$ be a regular Riemann Poisson manifold. If
$\al\in\Ga(Ker\pi)$, the following assertions are equivalent.

1) $\al$ is a basic 1-form.

2) $D\al=0$.

3) $\#(\al)$ is a foliate vector field.

4) $L_{\#(\al)}\pi=0.$

Furthermore, if $\al,\be\in \Om_b^1(P)$ then $<\al,\be>$ is a Casimir function.

\end{pr}
{\bf Preuve:}  By (16) and Proposition 1.1 we have $1)\ssi 2)$. 

$\#(\al)$ is a foliate vector field if and only if, for any
$\be\in\Ga((ker\pi)^\perp)$, for any $x\in P$ and $\ga_x\in Ker\pi_x$,
$\ga_x([\#(\al),\pi(\be)])=0.$

Let $\ga\in\Ga(Ker\pi)$. We have
\begin{eqnarray*}
\ga([\#(\al),\pi(\be)])&=&-d\ga(\#(\al),\pi(\be))-\pi(\be).\ga(\#(\al))\\
&=&-\pi(\be,i_{\#(\al)}d\ga)-\pi(\be,d[\ga(\#(\al))])\\
&=&-\pi(\be,L_{\#(\al)}\ga)=L_{\#(\al)}\pi(\be,\ga).\end{eqnarray*}
Now, by (12) and since obviously $L_{\#(\al)}\pi(\be,\ga)=0$ if $\ga,\be\in
\Ga(Ker\pi)$, we get $3)\ssi 4)$.

On other hand, we have also
$$\ga([\#(\al),\pi(\be)])=-d\ga(\#(\al),\pi(\be))-<D_\be\al,\ga>-<\al,D_\be\ga>.
$$
By Darboux theorem, for any $\ga_x\in Ker\pi_x$ there is a local Casimir
 function $f$ such that $d_xf=\ga_x$. This gives $2)\ssi 3).$

If $\la,\be$ are basic 1-forms, by (6) and 2) we have that $<\al,\be>$ is a
Casimir function.
$\Box$.

\begin{th} Let $(P,\pi,<,>)$ be a regular Riemann Poisson manifold. The
symplectic foliation is a K\"ahler-Riemann foliation and $g$ is bundle-like
metric.\end{th} 

{\bf Preuve:} We have shown that the symplectic leaves are K\"ahler. We show now
that $g$ is bundle-like. Following B. Reinhart [6], the metric $g$ is said to
be  bundle-like if it has the following property: for any open set $U$ in $P$ and
for all vector fields $Y;Z$ that are foliate and perpendicular to the leaves,
the function $g(Y,Z)$ is a basic function on $U$.

In our case, the perpendicular foliate vector fields are 
$Y=\#(\al)$ and $Z=\#(\be)$ where $\al,\be$ are basic 1-forms. Furthermore
$g(Y,Z)=<\al,\be>$ which is a Casimir function by Proposition 1.2 and the
theorem follows.$\Box$

{\bf Remark.} Let $(P,\pi,<,>)$ be a regular Riemann Poisson manifold. If the Lie algebra of foliate vector fields is transitif on $P$ ( which is the case if the symplectic foliation is transversally parallelizable)  the symplectic leaves are symplectomorph.

Now we give the converse of  Theorem 2.1.

\begin{th} Let $(P, F,g)$ be a differentiable manifold endowed with a
Riemannian foliation ${ F}$ and a bundle-like metric $g$. We suppose that
there is $\om\in\Ga(\wedge^2 F)$ such that:

1) for any leaf $L$, the restriction of $\om$ to $L$ is symplectic and parallel
with respect to the levi-Civita connection associated with the restriction of $g$
to $L$;

2) for any local perpendicular foliate vector field $X$  and any couple
$(U,V)$ of local vector fields tangent to $F$,
$$L_X\om(U,V)=0.$$ Then there is a Poisson tensor $\pi$ on $P$ and Riemann
metric $<,>$ on $T^*P$ such that $(P,\pi,<,>)$ is a regular Riemann Poisson
manifold whose symplectic foliation is $F$.\end{th}

{\bf Preuve:} We have
$$TP=F\oplus F',\qquad T^*P=F^\circ\oplus F'^\circ$$where $F'$
is the distribution $g$-orthogonal to $F$ and $F^\circ=\{\al\in T^*P/\al(F)=0\}.$

Denote  $\#:T^*P\too TP$ the identification given by the Riemann metric $g$. 
 We have
$$\#(F^\circ)=F'\quad\mbox{and}\quad \#(F'^\circ)=F.$$

The leafwise symplectic form $\om$ can be considered us a 2-form on $M$ by setting $i_v\om=0$ for any $v\in F'$ and hence realizes un isomorphism $\om:F\too
F'^\circ$, $v\mapsto\om(v,.)$. Denote $\om^{-1}:F'^\circ\too F$ its inverse.

Now, define a bivector $\pi$ by
$$\pi(\al,\be)=\left\{
\begin{tabular}{lcl}
$\om(\om^{-1}(\al),\om^{-1}(\be))$&if&$\al,\be\in F'^\circ$\\
0&if&$\al$ or $\be$ belongs to $F^\circ$,\end{tabular}\right.$$

and a Riemann metric on $T^*P$ by
$$<\al,\be>=\left\{
\begin{tabular}{lcl}
$g(\om^{-1}(\al),\om^{-1}(\be))$&if&$\al,\be\in F'^\circ$\\
$g(\#(\al),\#(\be))$&if&$\al,\be\in F^\circ$\\
0&if&$\al\in F'^\circ $ and $\be\in F^\circ$ .\end{tabular}\right.$$

Now, we will compute the Levi-Civita contravariant connection $D$ associated with
$(\pi,<,>)$. Before, we do some remarks on the Poisson bracket $[\;,\;]_\pi.$

For $\al,\be,\ga\in F'^\circ$ and $v$ a tangent vector to $F'$ in some $x$, a
straightforward calculation gives
\begin{eqnarray*}
[\al,\be]_\pi(\om^{-1}(\ga))&=&-d\om[\om^{-1}(\al),\om^{-1}(\be),\om^{-1}(\ga)]
-\om([\om^{-1}(\al),\om^{-1}(\be)],\om^{-1}(\ga)).\\
  \;[\al,\be]_\pi(v)&=&L_X\pi(\al,\be)\end{eqnarray*}where $X$ is any vector
field such that $X_x=v$. If we choose $X$ a foliate perpendicular vector field,
it 's easy to see that
$$L_X\pi(\al,\be)=L_X\om(\om^{-1}(\al),\om^{-1}(\be))=0.$$$d\om$ also vanishes
and we get
\begin{eqnarray*}
\;[\al,\be]_\pi(\om^{-1}(\ga))&=&
-\om([\om^{-1}(\la),\om^{-1}(\be)],\om^{-1}(\ga)).\\
\;[\al,\be]_\pi(v)&=&0.\end{eqnarray*}
From this relation we can deduce that
$$D_\al\be=\pi^{-1}(\na^L_{\pi(\al)}\pi(\be))$$where $\na^L$ is the Levi-Civita
connection of the restriction of $g$ to a leaf $L$.

If $\al,\be\in F'^\circ$ and $\ga\in F^\circ$, we have
\begin{eqnarray*}
[\al,\ga]_\pi(\pi(\be))&=&L_{\pi(\al)}\ga	(\pi(\be))=d\ga(\pi(\al),\pi(\be))\\
&=&-\ga([\pi(\al),\pi(\be)])=0\end{eqnarray*}since $F$ is involutif.

A straightforward calculation gives
$$2<D_\al\be,\ga>=\left\{
\begin{tabular}{lcl}
$0$
&if&$\al\in F^\circ$, $\be\in F'^\circ$, $\ga\in F'^\circ$.\\
$L_{\pi(\be)}g(\#(\al),\#(\ga))$&if&
$\al\in F^\circ$, $\be\in F'^\circ$, $\ga\in F^\circ$.\\
0&if&$\al\in F^\circ$, $\be\in F^\circ$, $\ga\in F^\circ$.\\
$L_{\pi(\ga)}g(\#(\al),\#(\be))$&if&
$\al\in F^\circ$, $\be\in F^\circ$, $\ga\in F'^\circ$.
\end{tabular}\right.$$
As a characterization of a bundle-like metric, we have
$$
L_{\pi(\ga)}g(\#(\al),\#(\be))=0\quad\mbox{ for}\qquad  \al,\be\in F^\circ, \ga\in
F'^\circ.\eqno(17)$$ 
Now, it is easy to verify that $D\pi=0$ and the theorem follows.$\Box$
\section{The Poisson cohomology of a regular Riemann Poisson manifold}

The Poisson cohomology of a Poisson manifold $(P,\pi)$ is the cohomology of the chain complex
$({\cal X}^*(P),d_\pi)$ where, for $0\leq p\leq dimP$, ${\cal X}^p(P)$ is the 
$C^\infty(P,\reel)$-module of $p$-multivector fields and $d_\pi$ is given by
\begin{eqnarray*}
d_\pi Q(\al_0,\ldots,\al_p)&=&\sum_{j=0}^p(-1)^j
\pi(\al_j).Q(\al_0,\ldots,\hat\al_j,\ldots,\al_p)\\
&+&\sum_{i<j}(-1)^{i+j}
Q([\al_i,\al_j]_\pi,\al_0,\ldots,\hat\al_i,\ldots,\hat\al_j,\ldots,\al_p).\qquad(18)
\end{eqnarray*}
We denote $H_\pi^*(P)$ the spaces of cohomology.

The leafwise cohomology of a manifold $P$ endowed with a foliation defined by an involutif
distribution $F$ is the cohomology of the chain complex $(\Om_{\cal F}^*(P),d_{\cal F})$ where,
for $0\leq p\leq rangF$,
$\Om_{\cal F}^p(P)$ is the $C^\infty(P,\reel)$-module of $p$-forms on the vector bundle $F$ and
$d_{\cal F}$ is given by
\begin{eqnarray*}
d_{\cal F} \om(X_0,\ldots,X_p)&=&\sum_{j=0}^p(-1)^j
X_j.\om(X_0,\ldots,\hat X_j,\ldots,X_p)\\
&+&\sum_{i<j}(-1)^{i+j}
\om([X_i,X_j],X_0,\ldots,\hat X_i,\ldots,\hat X_j,\ldots,X_p).\qquad(19)
\end{eqnarray*}
We denote $H_{\cal F}^*(P)$ the spaces of cohomology.

Let $(P,\pi,<,>)$ be a  Poisson manifold. 

For $1\leq p\leq dimP$, we consider 
the subspace 
${\cal X}^p_0(P)\subset{\cal X}^p(P)$ of $p$-multivector fields $Q$ such that $$i_\al Q=0\qquad\mbox{ for 
any}\qquad \al\in Ker\pi.$$
It's easy to verify that $d_\pi({\cal X}^p_0)\subset{\cal X}^{p+1}_0$. The natural injection ${\cal X}^p_0\inj{\cal X}^p$ induces a linear map $H^*({\cal X}^p_0)\too H_\pi^*(P)$ which is not injectif in general.

Let $\pi: 
\Om_{\cal F}^p(P)\too{\cal X}^p_0(P)$ be the map given by
$$\pi(\om)(\al_1,\ldots,\al_p)=\om(\pi(\al_1),\ldots,\pi(\al_p)).$$
It is easy to verify that $\pi$ is an isomorphism and $\pi(d_{\cal F}\om)=d_\pi\pi(\om)$ and
hence $\pi$ induces an isomorphism
$$\pi^*:H^p_{\cal F}(P)\too 
H^p({\cal X}^p_0(P))\eqno(20)$$where
$H^p_{\cal F}(P)$ is the leafwise cohomology of the symplectic foliation.

Suppose now that $(P,\pi,<,>)$ is a regular Riemann Poisson manifold.

  Let  
${\cal X}^p_1(P)\subset{\cal X}^p(P)$ be the space of $p$-multivector fields $Q$ such that
$$Q(\al_1,\ldots,\al_p)=0\qquad\mbox{for any } \al_1,\ldots,\al_p\in (Ker\pi)^\perp.$$ We have
$${\cal X}^p(P)={\cal X}^p_0(P)\oplus{\cal X}^p_1(P),\qquad 
d_\pi({\cal X}^p_i)\subset{\cal X}^{p+1}_i, i=0,1
\eqno(21)$$
and hence, for $1\leq p\leq dimP$, 
$$H^p_\pi(P)=H^p({\cal X}^*_0(P))\oplus H^p({\cal X}^*_1(P)).\eqno(22)$$

The isomorphism (20) gives an injectif linear map
$$\pi^*:H^p_{\cal F}(P)\inj 
H_\pi^p(P).\eqno(23)$$

We consider now the space $\Om^p_b(P)$ of   basic  $p$-forms. Let $\om\in
\Om^p_b(P)$, we denote $\#(\om)$ the $p$-multivector field in ${\cal X}^p_1(P)$ given by
$$\#(\om)(\al_1,\ldots,\al_p)=\om(\#(\al_1),\ldots,\#(\al_p)).$$
\begin{Le} Let $(P,\pi,<,>)$ be a regular Riemann Poisson manifold.
 For any basic $p$-form $\om$,
$$d_\pi\#(\om)=0.\eqno(24)$$\end{Le}
{\bf Preuve:} It's obvious that $d_\pi\#(\om)(\al_0,\ldots,\al_p)=0$ if all the $\al_i$ belong
to $\Ga(\ker\pi)$. 
On other hand since $\om$ is basic and since $\#$ maps $\Ga(Ker\pi)^\perp$ on $\Ga(TS)$,
$d_\pi\#(\om)(\al_0,\ldots,\al_p)=0$ if there is $i\not=j$ such that $\al_i$ and $\al_j$ belong to
$\Ga((Ker\pi)^\perp)$.

 Now, we suppose that $\al_0\in \Ga((Ker\pi)^\perp)$ and, for 
$1\leq j\leq p$, $\al_j\in
\Ga(\ker\pi)$. We have
\begin{eqnarray*}
d_\pi\#(\om)(\al_0,\ldots,\al_p)&=&\pi(\al_0).\om(\#(\al_1),\ldots,\#(\al_p))\\
&+&\sum_{j=1}^p(-1)^j\om(\#([\al_0,\al_j]_\pi),\#(\al_1),\ldots,\hat\al_j,\ldots,\#(\al_p)).
\end{eqnarray*}
Now we will compute the perpendicular component of $\#([\al_0,\al_j]_\pi)$. Let $\ga\in
\Ga((Ker\pi)$.
\begin{eqnarray*}
\ga[\#([\al_0,\al_j]_\pi)]&=&\ga[\#(L_{\pi(\al_0)}\al_j)]=<L_{\pi(\al_0)}\al_j,\ga>\\
&=&-<\al_j,L_{\pi(\al_0)}\ga>+\pi(\al_0).<\al_j,\ga>\qquad \mbox{by}\quad (17)\\
&=&-<\al_j,i_{\pi(\al_0)}d\ga>+\pi(\al_0).<\al_j,\ga>\\
&=&-d\ga(\pi(\al_0),\#(\al_j))+\pi(\al_0).<\al_j,\ga>\\
&=&\ga([\pi(\al_0),\#(\al_j)]).\end{eqnarray*}
So $[\pi(\al_0),\#(\al_j)]
$ is the perpendicular component of $\#([\al_0,\al_j]_\pi)$ and
$$d_\pi\#(\om)(\al_0,\ldots,\al_p)=L_{\pi(\al_0)}\om(\#(\al_1),\ldots,\#(\al_p))=0$$
since $\om$ is basic.$\Box$

We get a linear map $\#:\Om^p_b(P)\too H^p({\cal X}^p_1(P))$ which is obviously injectif. We have shown the following theorem.
\begin{th} Let $(P,\pi,<,>)$ be a regular Riemann Poisson manifold. For $1\leq p\leq dimP$, there is an injectif linear map $$
\Om^p_b(P)\oplus H^p_{\cal F} (P)\inj H^p_\pi(P ).\eqno(25)$$
In particular, we have
$$
\Om^1_b(P)\oplus H^1_{\cal F}(P)\simeq H^1_\pi(P ).\eqno(26)$$
\end{th}

{\bf Remarks.} 

 Let $(P,\pi,<,>)$ be a regular Riemann Poisson manifold. On can verify that the injection  $H^1_{\cal F}(P)\inj H^1_\pi(P)$ maps
the modular class of the symplectic foliation which is an obstruction lying in $H^1_{\cal F}(P)$ to the existence of an invariant transverse volume form to $2mod(P)$ where $mod(P)$ is the modular class of the Poisson manifold which is an obstruction lying $H^1_{\pi}(P)$ to the existence of a volume form invariant with respect to Hamiltonian flows. Since the modular class vanishes for a Riemannian foliation, we get that $(P,\pi)$ is unimodular and we have another proof of the result given in [2].

As a final we give some examples:

1) Any Poisson manifold such that the Poisson tensor is parallel with respect to the Levi-Civita connection for a Riemannian metric is a regular Riemann Poisson manifold (  see [1]).

2) Let $(M,\om,g)$ be a K\"ahler manifold endowed with a symplectic, isometric and free action of the circle $S^1$. Let $Z$ denote de fundamental vector field of this action and suppose that there is a vector field $Y$ such that $Z$ and
$Y$ commute and  are orthogonal. The 1-form $i_Z\om$ project on the space of orbits $M/S^1$ and define a Riemannian foliation which verify the hypothesis of Theorem 2.2.\bigskip\eject

{\bf References}

\bigskip

[1] M. Boucetta,  Compatibilit\'e des structures pseudo-riemanniennes et des
structures de Poisson, C. R. Acad. Sci. Paris, t. 333, S\'erie I, p. 763-768,
2001.

[2] M. Boucetta, Poisson manifolds with compatible pseudo-metric and pseudo-Riemannian Lie algebras, Preprint math.DG/0206102.

[3] C. Deninger and W. Singhof, Real Polarizable Hodge Strutures Arising from foliations, Annals of Global Analysis and Geometry 21, 377-399, 2002.

[4]  R. L. Fernandes, Connections in Poisson Geometry 1: holonomy and
invariants, J. Diff. Geom. 54 , p. 303-366, 2000.

[5] P. Molino, Riemannian foliations, Progress in Mathematics Vol. 73, Birkh\"auser Boston Inc. 1988.

[6] B. Reinhart, Foliated manifolds with bundle-like metrics, Ann. Of Math. 81  119-132, 1959.

[7] I. Vaisman, Lectures on the Geometry of Poisson Manifolds, Progress in
Mathematics, vol. 118, Birkh\"auser, Berlin, 1994.

[8] A. Weinstein, The Modular Automorphism Group of a Poisson Manifold, J.
Geom. Phys. 23, p. 379-394, 1997.

[9] A. Weinstein, The Local Structure of Poisson Manifolds, J. Differential
Geometry 18, p. 523-557, 1983.

\bigskip

{\bf Mohamed BOUCETTA

Facult\'e des sciences et techniques, BP 549, Marrakech, Maroc

Email: boucetta@fstg-marrakech.ac.ma}

\end{document}